\documentclass{article}
\usepackage{amsthm,amsfonts, amsbsy, amssymb,amsmath,graphicx}
\usepackage{graphics}
\newtheorem{thm}{Theorem}
\newtheorem{lm}{Lemma}

\newtheorem{re}{Remark}

\author{Vassily Olegovich Manturov}

\title{Khovanov homology for virtual knots with arbitrary coefficients}

\begin{document}

\maketitle

\begin{abstract}We construct explicitly the Khovanov homology theory for virtual links
with arbitrary coefficients by using the twisted coefficients
method. This method also works for constructing Khovanov homology
for ``non-oriented virtual knots'' in the sense of \cite{Vir}, in
particular, for knots in ${\bf R}P^{3}$.
\end{abstract}

\newcommand{\cC}{{\cal C}}

\begin{figure}
\centering\includegraphics[width=360pt,height=135pt]{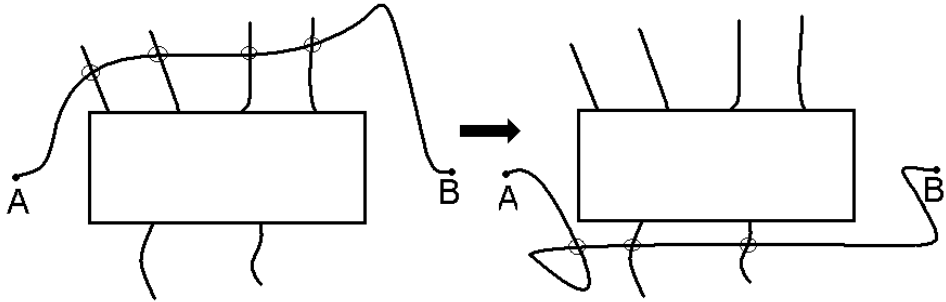}
\caption{The detour move} \label{detour}
\end{figure}

Virtual knots were introduced in mid-nineties by Lou Kauffman, see
\cite{KaV}. By a virtual diagram we mean a four-valent graph on
the plane endowed with a special structure: each crossing is
either said to be classical (in this case one indicates which pair
of opposite edges at this crossing forms an overcrossing; the
remaining two edges form an undercrossing) or virtual (in this
case, we do not specify any addtional structure; virtual crossings
are just marked by a circle). Two virtual diagrams are called {\em
equivalent} if one of them can be obtained from the other by a
finite sequence of generalized Reidemeister moves and planar
isotopies. Recall that the list of generalized Reidemeister moves
consists of classical Reidemeister moves (see,e.g., \cite{Man})
and the detour move. The latter means that if wa have a purely
virtual arc containing only virtual crossings, we may remove it
and restore in any other place of the plane, see Fig.
\ref{detour}.

One of the most important achievements in modern knot theory is
the Khovanov categorification of the Jones polynomial proposed in
\cite{Kh}. The main idea is to replace (Laurent) polynomials by
homology groups of bigraded vector spaces. Thus, with each link
diagram, one associates a bigraded algebraic complex; all
cohomology groups of this complex are knot invariants, whence its
graded Euler characteristic coincides with the Jones polynomial
(in a certain normalization).

Later on, we deal with bigraded complexes of the type
$\cC=\oplus_{i,j}\cC^{i,j}$ with {\em height} $i$ and {\em
grading} $j$. The differential $\partial$ in such a complex
preserves the grading and increases the height by one.

For a graded linear space $A$ define its {\em graded dimension}
$qdim$ by $qdim A=\sum_{i} q^{i}\cdot dim A^{i}$, where $A^{i}$ is
the subspace of  $A$ generated by all vectors of grading $i$.
Chain spaces  $\cC^{j}$ of $\cC$ of a given height $j$ can be
treated as such graded spaces.

It would be more conveninent to speak about {\em cohomology}, but
we should rather say {\em Khovanov homology} instead and say {\em
chains, cycles, boundaries} rather than {\em cochains, cocycles,
coboundaries}. For such complexes there are well defined degree
shift and height shift operations: $\cC\mapsto \cC[k,l]$; where
$(\cC[k,l])^{i,j}=\cC^{i-k,j-l}$. By the {\em graded Euler
chararcteristic} of a complex $\cC^{ij}$ we mean the alternating
sum of graded dimensions of the chain spaces, or, equivalently,
that of homology groups. For the chain spaces we have:
$\sum_{i}(-1)^{i}qdim {\cC}^{i}=\sum_{i,j} (-1)^{i}q^{j} dim
{\cC}^{i,j}$.

A beautiful property of the Khovanov homology is its {\em
functoriality} meaning the following. Each cobordism in  ${\bf
R}^{3}\times I$ between two links $K_{1}\subset {\bf
R}^{3}\times\{0\}$ and $K_{2}\subset {\bf R}^{3}\times \{1\}$,
generates a natural mapping between the Khovanov homology
$Kh(K_{1})\to Kh(K_{2})$ which is invariant under isotopy. Here
the Khovanov homology is built of bricks, which are circles on the
plane, and the complex cobordism category can be viewed as an
algebraic counterpart of these circles where maps are represented
by cobordisms.

The functoriality of the Khovanov homology was proved by Jacobsson
in the algebraic setup and then by Bar-Natan in the topological
setup, see \cite{BN2,Jac}.

Another definition of the Khovanov homology is {\em local}:
instead of considering states and counting the number of circles
in each of them, one can construct a complex by using Matrix
factorizations at every vertex, as in \cite{KhR1}, and then
contracting them into the whole complex. The homology proposed in
\cite{KhR1} and \cite{KhR2} give a categorification of the HOMFLY
polynomial, the Jones 1-variable polynomial being a partial case
of it.

This way is, in some sense, easier to check the invariance
(because the proof becomes local), but it is much more complicated
for explicit calculations. However, it is said on page 32 of
\cite{KhR1} ``The network does not even have to be planar, and
does not need to be embedded anywhere. In our paper, however, all
such diagrams are going to be planar''. This  means that the
construction of \cite{KhR1} [not only for the usual
$sl(2)$-homology, but in the general case] generalizes for the
case of virtual knots. Indeed, all their proofs are local, and the
construction of embedded graph does not feel nugatory handles of a
$2$-surface, thus being invariant under (de)stabilizations.

However, this construction is too much implicit. In this paper, we
present an explicit construction of this homology together with
several generalizations of it.

It turns out, for example, that the Khovanov homology for virtual
knots constructed in this way is invariant under virtualizations
thus leading to a homology theory of ``twisted virtual knots''.

One can also establish the analogues of some other theories: Lee's
theory, Wehrli's spanning tree expansion, etc.

In the previous papers \cite{Khoz2,Kho}, the author constructed
the Khovanov theory for virtual knots. More precisely, the
Khovanov complex was well defined for all virtual link diagrams
only over ${\bf Z}_{2}$; the Khovanov homology with arbitrary
coefficients was defined only in the case of virtual knots
corresponding to so-called {\em orientable atoms}: {\bf atom plays
the key role in defining the Khovanov homology for virtual links,
with its non-orientability being the main obstruction}. For
non-orientable atoms, we presented two geometrical construction
transforming them to orientable atoms, and then proving that the
Khovanov homology of the target knot (atom) is an invariant of the
source one.

In the present note, we will, for arbitrary coefficient ring,
construct a differential with $\partial^{2}=0$ explicitly. The key
ideas are to use the {\bf twisted coefficient} methods: one should
change the basis of the Hopf algebra (representing the Khovanov
homology of the unknot) while passing from one classical crossing
to another crossings and to use the exterior product of tensor
spaces instead of the usual symmetric product.

Thus, with each virtual link diagram we associate a bigraded
complex with homology being invariant under generalized
Reidemeister moves. We wish to point out the main properties of
this construction.

\begin{enumerate}

\item The complex is constructed by using {\em atoms}; this it is
invariant under {\em virtualization}, the natural transformation
of virtual links preserving the Jones polynomial.

\item Since there exists a map from twisted virtual knots to usual
virtual knots modulo virtualization, this approach works for
twisted virtual knots as well.

\item In the ${\bf Z}_{2}$-case the complex coincides with the
${\bf Z}_{2}$-complex first proposed in \cite{Khoz2}.

\item For orientable atoms this complex has the same homology as
the complex proposed in \cite{Kho}.

\item The invariance proof is local and repeats that for the
classical case, see, e.g., \cite{BN}, the main difficulty being
the definition of the complex.

\item A partial case of this theory is theory of knots in ${\bf
R}P^{3}$; recall that the Kauffman bracket version of the Jones
polynomial for such knots was proposed in \cite{Dro}.

\end{enumerate}

\section{The Kauffman bracket and the Jones polynomial. Atoms. Twisted virtual knots}

Consider an oriented virtual diagram  $L$ and a diagram $|L|$
obtained from $L$ by forgetting the oritentation. Let us smooth
classical crossings of  $|L|$ according to the following rule.
Each classical crossing can be {\em smoothed} in one of two ways,
the positive way $A$ and the negative way $B$.

A {\em state} is a choice of smoothing for all classical crossings
of a diagram. Each state generates a set of curves on the plane
having only virtual crossings. In other words, we get some virtual
unlink. Suppose the diagram $L$ has $n$ classical crossings.
Enumerate them arbitrarily. We get  $2^{n}$ states which can be
encoded by vertices of the $n$-dimensional discrete cube
$\{0,1\}^{n}$, where $0$ and $1$ correspond to $A$-smoothing and
$B$-smoothing, respectively. Call this cube the {\em state cube}
of the diagram. Two states are {\em adjacent} if they differ in a
unique coordinate. Any two adjacent vertices are connected by an
edge of the cube. Orient all edges of the cube as the coordinate
increases, i.e., from  $A$-type smoothing to $B$-type smoothing.
By {\em height} of the vertex we mean the number of $B$-type
smoothing of the corresponding state.

For each state $s$, denote by $\alpha(s)$ the number of
$A$-smoothings, let  $\beta(s)=n-\alpha(s)$ be the number of
$B$-smoothings and let $\gamma(s)$ be the number of componentsof
the resulting unlink. The Jones-Kauffman polynomial is then
defined as:

\begin{equation}X(L)=(-a)^{-3w(L)}\sum_{s}a^{\alpha(s)-\beta(s)}(-a^{2}-a^{-2})^{\gamma(s)-1},
\label{18n1}\end{equation} where $w(L)$ is the writhe of the
oriented diagram $L$ (i.\,e., the difference between the number of
positive crossings and the number of negative crossings).

The unnormalized version
$\sum_{s}a^{\alpha(s)-\beta(s)}(-a^{2}-a^{-2})^{\gamma(s)-1}$ is
called  the {\em Kauffmab bracket}; both the Kauffman bracket and
the Jones-Kauffman polynomial are (Laurent) polynomials in one
variable.  The Jones-Kauffman polynomial is invariant under
(generalized) Reidemeister moves.

\newcommand{\Jj}{\hat J}

After the variable change $a=\sqrt{(-q^{-1})}$, we get instead of
the Jones polynomial $V$ its modified versions $J$ and $\Jj$. They
differ by a normalization on the unknot: we have $J=1$ on the
unknot whence $\Jj=1$ for the empty link with no components.
Herewith $J=\frac{\Jj}{(q+q^{-1})}$. In more details,  $\Jj$ goes
as follows. Let $L$ be an oriented virtual diagram and let $|L|$
be the corresponding unoriented virtual diagram; denote by $n_{+}$
and $n_{-}$ the number of positive (resp., negative) classical
crossings of $L$; let  $n$ be the total number of classical
crossings: $n=n_{+}+n_{-}$. Set

\begin{equation}\Jj(L)=(-1)^{n_{-}}q^{n_{+}-2 n_{-}}\langle
L\rangle,\label{nenorm}\end{equation} where $\langle L\rangle$ is
the modified Kauffman bracket defined axiomatically as $\langle
\bigcirc\rangle=(q+q^{-1})$, $\langle L \rangle=\langle L_{A}
\rangle-q\langle L_{B} \rangle$, $\langle L\sqcup
\bigcirc\rangle=(q+q^{-1})\cdot \langle L\rangle$.

Later on, we deal with the polynomial $\Jj$ calling it {\em the
Jones polynomial}.

The polynomial $\Jj$ has an easy description in terms of the state
cube. Taking away the normalization factors
$(-1)^{n_{-}}q^{n_{+}-2 n_{-}}$ we get the (slightly modified)
Kauffman bracket
$\sum_{s}(-q)^{\beta(s)}((q+q^{-1})^{\gamma(s)})$. That is, we
take the sum over all vertices of the cube of  $(-q)$ to the power
equal to the {\em height of the vertex}, all multiplied by
 $(q+q^{-1})$ to the number of circles corresponding to the vertex of the
cube. Thus, on the level of polynomials, with each circle we
associate  $(q+q^{-1})$ and then we multiply these polynomials,
shift them by multiplying by  $\pm q^{k}$ and collect as summands.

Thus, the Jones polynomial can be restored from the information
about the {\em number of circles} at the states of the cube. If we
also take into account the way {\bf how these circles bifurcate}
while passing from a state to the adjacent state, we get the
Khovanov complex.

It turns out that all this information needed for the Jones
polynomial is contained in the atom representing its planar
diagram.

An {\em atom} is a pair $(M,\Gamma)$ : closed two-manifold
 $M$ and a four-valent graph (frame) $\Gamma\subset M$ dividing $M$ into
 black and white cells together with a checkerboard coloring of these cells. Atoms
 are considered up to natural combinatorial equivalence: diffeomorphisms
 of the underlying manifolds preserving the frame and the cell coloring.

Associate with a given virtual link diagram  $L$ the atom  $V(L)$
to be constructed as follows. The vertices of  $V(L)$ are in
one-to one correspondence with {\em classical} crossings of  $L$.
These classical crossings are connected by {\em branches of the
diagram}, which might intersect on the projection plane. For each
classical crossing we have four emanating branches. With these
edges we associate four edges of the atom frame to connect the
corresponding vertices. The rule for attaching black (resp.,
white) cells is defined by the planar diagram $L$. Namely, let $X$
be a classical crossing of $L$. Enumerate the four emanating edges
by $x_{1},x_{2},x_{3},x_{4}$ in the clockwise direction in such a
way that the branches  $x_{1}$ and $x_{3}$ form an undercrossing
whence  $x_{2}$ and $x_{4}$ form an overcrossing. Then for
attaching the black cells we chose those pairs of half-edges of
the atom corresponding to $(x_{1},x_{2})$ and $(x_{3},x_{4})$.

Evidently, the whole information about the number of circles in
states of the diagram can be extracted from the corresponding
atom. In other words, the state cube can be completely restored
from the atom.

By {\em virtualization} we mean the local transfomation in the
neighbourhood of a classical crossing shown in Fig \ref{twist}.
Note that the virtualization does not change the state cube; it
does not change the atom, either. Moreover, two diagrams
corresonding to the same atoms, can be obtained from each other by
a sequence of detours and virtualizations.

\begin{figure}
\centering\includegraphics[width=180pt, height=100pt]{twist.bmp}
\caption{Two variants of virtualization} \label{twist}
\end{figure}

We shall construct the Khovanov complex starting from the atom
corresponding to a given diagram. Thus the homology will be
automatically invariant under virtualization. This supports the
conjecture stated in \cite{FKM} that if two calssical crossings
are equivalent by a chain of generalized Reidemeister moves,
detours and virtualizations, then the corresponding links are
isotopic in the usual sense.

Twisted virtual knots \cite{Vir} are close relatives of virtual
knots. They are represented by knots in oriented thickenings of
not necessarily orientable surfaces modulo
stabilization/destabilization.

Namely, having a non-orientable surface $S$, one can construct the
canonical oriented $I$-bundle over it, which is a $3$-manifold
$S{\tilde \times} I$ with boundary.

A nice example of such a thickened surface is ${\bf R}P^{2}{\tilde
\times }I$, which is homeomorphic to ${\bf R}P^3$. Thus,
constructing a Khovanov homology for such knots, we shall get a
Khovanov homology theory for knots in ${\bf R}P^{3}$.

Given a surface $S$, one can project any knot (link) from
$S{\tilde\times}I$ to $S$. In the general position, one gets a
$4$-graph. In order to restore the link, one should indicate for
which crossings how the two branches behave. In the orientable
case, one just indicates, which branch should be over, and which
branch should be under. However, in the non-orientable case this
indication is relative: while walking along a non-orienting
circuit, the direction upwards changes to the direction downwards.

However, such surfaces perfectly agree with atoms. Indeed, fix
(once forever) an orientation of $S{\tilde\times} I$. Now, for
link diagram in $S$, we already have a frame of the atom: a
$4$-valent graph with $A$-structure.

Now, the way for attaching black cells is the following: for a
vertex $X$, take two emanating edges $a$ and $b$. In the
neighbourhood of $X$, denote the small vector going from the edge
$a$ to the edge $b$. Now, if the orientation $(a,b,c)$ coincides
with the orientation of our $3$-manifold, then the angle generated
by the vectors $a$ and $b$ is black, as well as the angle opposite
to it; the remaining to angles are white.

Alternatively, if the orientations disagree, the angle between $a$
and $b$ is white.

Note that this choice does not depend on the ordering of $a$ and
$b$, nor on their directions.

\begin{figure}
\centering\includegraphics[width=200pt]{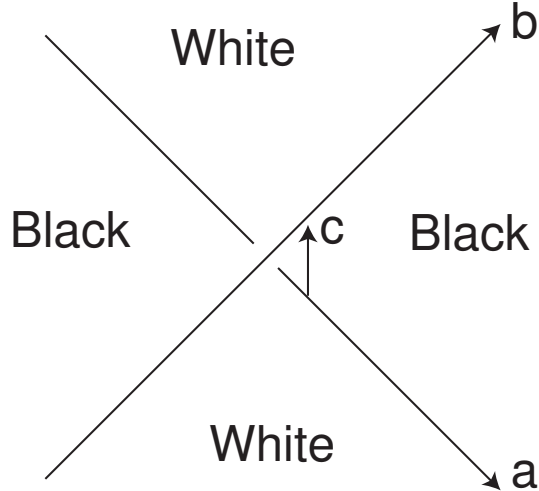}
\caption{Constructing an atom from a diagram} \label{skz}
\end{figure}

This leads to the following
\begin{thm}
There is a well-defined map from twisted virtual knots to virtual
knots modulo virtualization.\label{vzt}
\end{thm}

Knots in such surfaces were considered by Asaeda, Przytycki and
Sikora in \cite{APS}, and Viro \cite{Vir}. In \cite{APS} a
Khovanov homology theory for such surfaces was constructed by
using an additional topological information coming from surfaces.

From Theorem \ref{vzt} it follows immediately that the Khovanov
homology constructed below generalizes for twisted virtual knots.

\section{Definition of the Khovanov homology}

Our aim is to define a homology theory for virtual knots
(orientable and unorientable) over arbitrary ring in such a way
that:

\begin{enumerate}

\item The homology we are defining is invariant under Reidemesiter
moves

\item For the case of orientable atoms (alternatible knots) this
homology theory coincides with the one constructed in \cite{Kho}.

\item The ${\bf Z}_2$-restriction of this homology theory
coincides with the theory constructed in \cite{Khoz2}.

\end{enumerate}

In \cite{Kho} we constructed a homology theory for orientable
virtual knots with arbitrary coefficients. The main obstruction to
extend this theory over unorientable atoms is the possibility of
$1\to 1$ bifurcation on an edge of the cube.

If no such bifurcations occur, we may construct the Khovanov cube
by using the standard differentials, the multiplication $m$ (for
$2\to 1$-bifurcations) and the comultiplication $\Delta$ (for
$1\to 2$-bifurcations).

The situation with $1\to 1$ bifurcation makes the problem more
complicated. Indeed, if we wish to construct a grading-preserving
theory without introducing any new grading, this differential
should be identically equal to zero because of the dimension
reasons (there should be a mapping from $V$ to $V$ that lowers the
grading by one).

If we set this differential to be equal to zero with all other
differentials ($m$ and $\Delta$) defined in the standard way, we
get a straightforward generalization for the ${\bf Z}_{2}$ case.

The case of ${\bf Z}$-coefficients is more delicate: a
straightforward generalization makes some $2$-faces of the cube
commute and some other ones anti-commute.

{\bf Notational agreement}. Given several vector spaces
$V_{1},\dots, V_{n}$. We shall distinguish between two types of
tensor products, the {\bf ordered one} and the {\bf unordered
one}. For any permutation of indices $\sigma_{1},\dots,
\sigma_{n}$, we may consider the tensor product
$V_{\sigma_1}\otimes \dots V_{\sigma_{n}}$.

In the unordered case, we identify all these tensor products by
$x_{\sigma_{j}}\otimes \dots \otimes
x_{\sigma_{n}}=x_{1}\otimes\dots\otimes x_{n}$, where $x_{k}\in
V_{k}$. In the second case, we set $x_{\sigma_{j}}\otimes \dots
\otimes x_{\sigma_{n}}=sign(\sigma) x_{1}\otimes \dots\otimes
x_{n}$.

We shall denote such a tensor product of elements by the usual
tensor product sign $x_{1}\otimes x_{2}\otimes \dots \otimes
x_{n}$ in the first case and by using the exterior product sign
$x_{1}\wedge x_{2}\wedge \dots \wedge x_{n}$.

\begin{re}
To avoid confusion, note that, writing $X\wedge X$, we always
assume that the first $X$ and the second $X$ belong to different
spaces; thus $X\wedge X$ is not zero; one should rather write
$X_{1}\wedge X_{2}=-X_{2}\wedge X_{1}$, where $X_{i}$ here means
the element $X$ lying in the space $V_{i}$.
\end{re}

To handle it and to make the whole cube anti-coomutative, we have
to add two ingredients, sensitive to orientability of the atom:

\begin{enumerate}

\item With each circle of any state, we associate a vector space
of graded dimension $1+q^{2}$, however, without a fixed basis.
More precisely, we fix one element of the basis, denoted by $1$,
of grading $0$. The other element of the basis is defined up to
$\pm$ sign.

Given an orientation $o$ of the circle $C$ , we associate with it
the graded vector space generated by $1$ and $X_{C,o}$. If we
change the orientation to the opposite one $-o$, we set
$X_{C,-o}=-X_{C,o}$.

\item Given a state $s$ of the cube with $k$ circles $C_{1},\dots,
C_{k}$ in it. With this state, we associate the ordered tensor
product $V^{\otimes k}$ to it, the basis elements of this space
being products $(p^{1})_{C_{a_1}}\wedge
(p^{2})_{C_{a_2}}\wedge\dots\wedge (p^{k})_{C_{a_k}}$, where each
$(p^{i})$ is an element of the basis of $V_{C_{a_k}}$; here
$p_{C_i}\wedge q_{C_j}=-q_{C_j} \wedge p_{C_i}$.
\end{enumerate}

\begin{figure}
\centering\includegraphics[width=300pt]{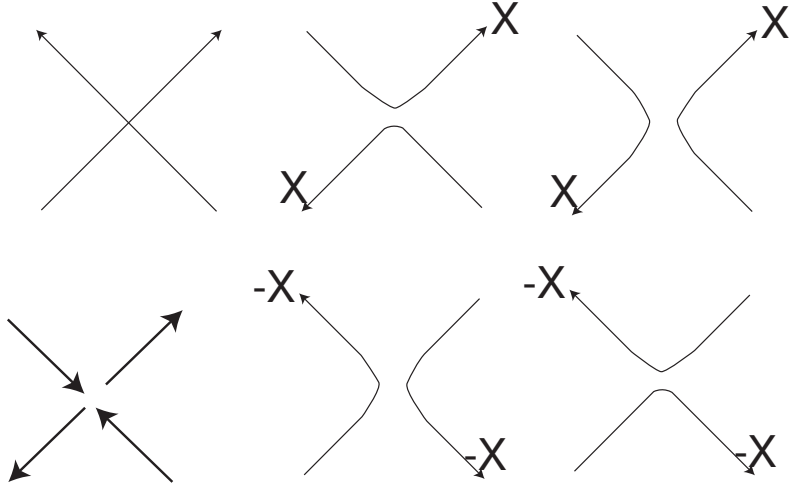}
\caption{Determining the basis for a given vertex} \label{figj}
\end{figure}

Now, the differentials are defined with respect to the
orientations of the circles at vertices and local ordering of the
components, as follows. For each vertex $v$, we fix the
orientations of the circles incident to this vertex according to
the orientation of the upper-right outgoing edge, see Fig
\ref{figj}.

That is, the orientation of these circles agrees with the
orientation of upper-right and lower-right edges on the knot
diagram and disagrees with that for the upper-left and lower-left
edges: {\bf We orient half-edges incident to a crossing as shown
in the lower-left picture of Fig. \ref{figj}; we say that the
orientation of a state circle incident to it is positively
oriented if it agrees with the local orientations indicated in the
lower-left picture; otherwise, the circle is negatively oriented.}

This orientation is well-defined unless the edge corresponding to
this vertex transforms one circle to one circle. {\bf For such
situations, we do not need to define $X$ at this vertex; we just
set the partial differential corresponding to this edge to be
identically zero}.

Assume we have a $1\to 2$ or $2\to 1$-bifurcation at a vertex $X$.

If we deal with two vertices adjacent to a vertex from opposite
sides, we enumerate them so that the upper [resp., left] one is
denoted as the first one, and the lower [resp., right] one is the
second one.

{\bf Agreement.} Later on, while defining $m$ and $\Delta$ we
assume that the circles we deal with are in the very first
position in our ordered tensor product. Otherwise, we permute them
to get the desired expression, possibly, with minus sign.

Now, we define $\Delta$ and $m$ locally with respect to the
prescribed choice of generators in $V$'s and the prescribed
ordering.

$\Delta (1)=1_{1}\wedge X_{2}+X_{1}\wedge 1_{2};
\Delta(X)=X_{1}\wedge X_{2}$ and

\begin{figure}
\centering\includegraphics[width=300pt]{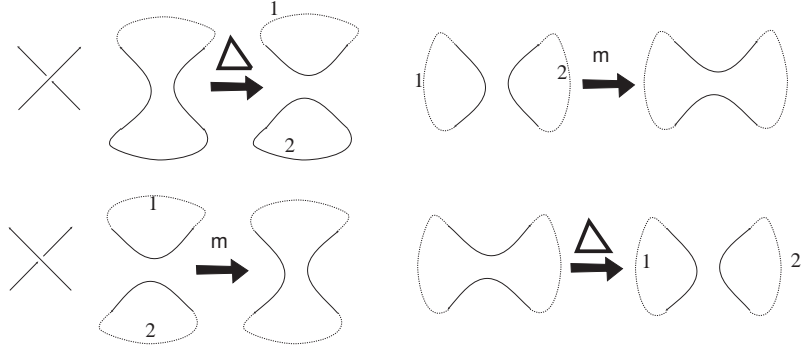}
\caption{Definition of $m$ and $\Delta$} \label{defmD}
\end{figure}

$m(1_{1}\wedge 1_{2})=1;m(X_{1}\wedge 1_{2})=m(1_{1}\wedge
X_{2})=X; m(X_{1}\wedge X_{2})=0$, see Fig. \ref{defmD}.

Thus, if we wish to comultiply the second factor $X_{2}$ in
$X_{1}\wedge X_{2}$, we get $X_{1}\wedge X_{2}=-X_{2}\wedge
X_{1}\to - X_{2}\wedge X_{3}\wedge X_{1}=-X_{1}\wedge X_{2}\wedge
X_{3}$, where $X_{3}$ belongs to the new (third) component.

\section{Formulation and proof of the main theorem}

Denote the obtained complex for a virtual diagram $K$, by $[[K]]$.

\begin{thm}
The complex $[[K]]$ described above, is well-defined for every
virtual diagram $K$.
\end{thm}

We have to prove that the complex we obtain is well-defined.

Then, as in \cite{Kho}, we normalize it by setting ${\cal
C}:=[[L]][-n_{-}]\{n_{+}-2n_{-}\}$. The obtained complex (that
differs from the unnormalized one just by a grading and degree
shift) will be well-defined as well. We will show that  its
homology is invariant under Reidemeister moves, this homology
coincides with the one constructed in \cite{Kho} in the orientable
case.

When considered modulo $2$, the definition agrees with the one
given in \cite{Khoz2} verbatim.

We first prove two lemmas establishing some properties of the
complex and making the further arguments easier.

Given a virtual diagram $K$ and a classical crossing $V$ of it.
Consider the diagram $K'$ obtained from $K$ by virtualizing $V$.
There is a natural one-to-one correspondence between classical
crossings of $K$ and $K'$. For each classical crossing $U$ of $K$,
we denote the space corresponding to it by $C_{U}$; denote the
corresponding space for $K'$ by $C_{U'}$.

\begin{lm}
Let $K,K'$ be two knot diagrams obtained one from another by a
virtualization. Then there is a grading-preserving chain map
$C(K)\to C(K')$ that maps $C_{U_K}$ isomorphically to $C_{U_{K'}}$
and agrees with the local differentials.

In particular, if $C(K)$ is a well-defined complex, then so is
$C(K')$ and their homology groups are isomorphic. \label{lem1}
\end{lm}

\begin{proof}
Suppose a diagram $K'$ is obtained from a diagram $K$ by
virtualization at a crossing $U$. Consider the corresponding cubes
$C(K)$ and $C(K')$. Obviously, their differentials agree for all
edges corresponding to all crossing except $U$. We assume the
differentials corresponding to $U$ divide our cubes into ``top
layer'' and the ``bottom layer'', as shown in Fig. \ref{lls}.

\begin{figure}
\centering\includegraphics[width=300pt]{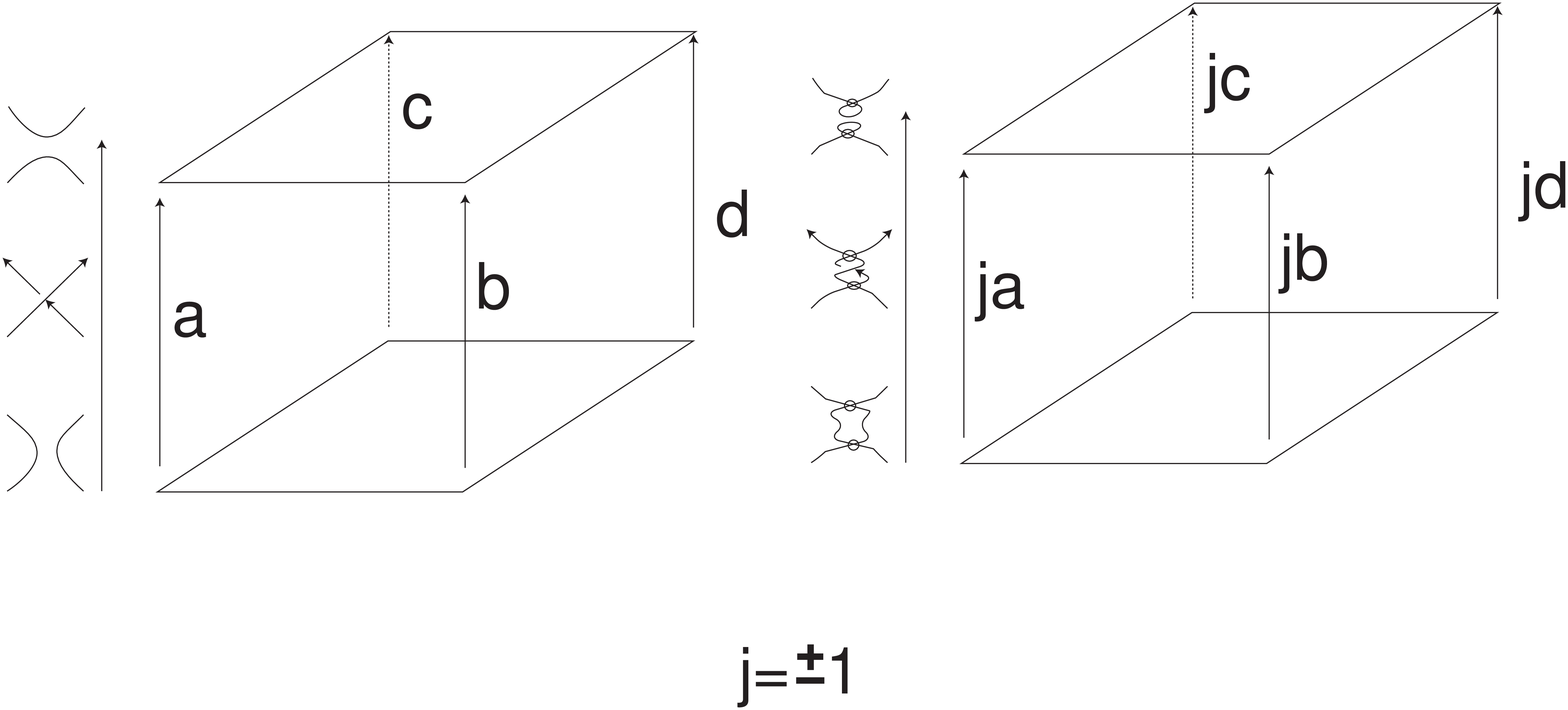}
\caption{Behaviour of the cube under virtualization} \label{lls}
\end{figure}

Now, the remaining differentials  differ possibly by signs on
edges corresponding to the crossing $U$. Our goal is to show that
they either all agree or all differ by $-1$ sign, as shown in Fig.
\ref{lls}.

Indeed, the bases at all crossings but $U$ agree for $K$ and $K'$.
This leads to the identification of chains of the corresponding
complexes. Under this isomorphism, we see that for every circle
incident to the crossing $U$, $X_{U,K}=-X_{U,K}$.

If we dealt with the usual tensor product case regardless the
circle ordering, the transformation $X\to -X$ would leave $m$
invariant and change $\Delta$ to $-\Delta$.

But we shall also take into account the circle ordering at a
vertex.

First assume $X$ is positive. Then all mappings $m$ corresponding
to $X$ represent a multiplication of elements corresponding to two
circles, a left one and a right one. After a virtualization, the
former left circle becomes the right one, and the former right one
goes to the left, see Fig. \ref{LR}.

\begin{figure}
\centering\includegraphics[width=300pt]{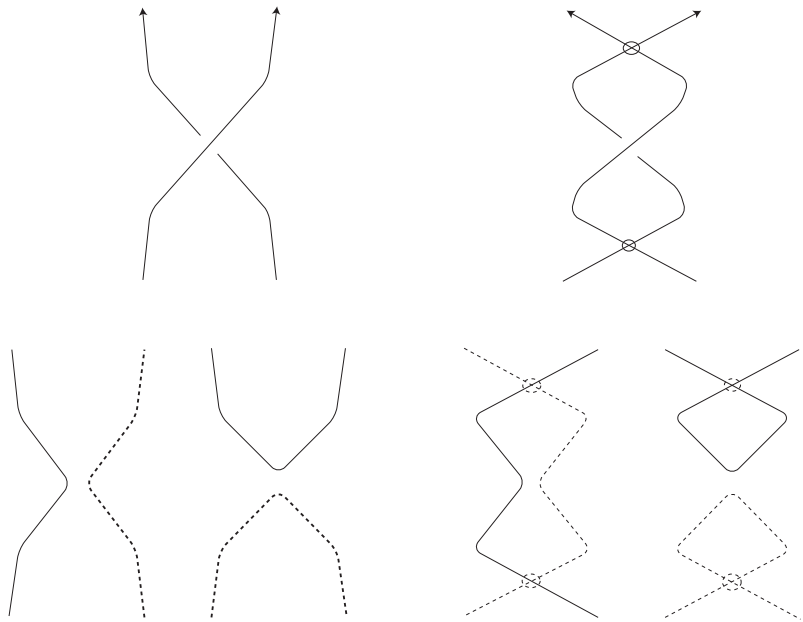} \caption{A
virtualization} \label{LR}
\end{figure}

Globally, thus, it changes the sign of $m$ in the fixed basis. For
the case of $\Delta$, we deal with one circle which is transformed
to two circles, an upper one and a lower one; the relation
``upper-lower'' is preserved by virtualization. The first
component above is represented by a solid line; the second
component is dashed.

Summing up, we see that a virtualization at a positive crossing
changes all signs of local differential corresponding to this
vertex.

Now, divide the cubes $C(K)$ and $C(K')$ into two parts with
respect to the smoothing of $U$, the lower one and the upper one.
Define the mapping $C(K)\to C(K')$ to be identity for all elements
from the lower cube and minus identity from all elements of the
upper cube.

Evidently, this mapping agrees with local differentials, and, if
the initial cube were commutative, this mapping would provide a
homology isomorphism.

A similar argument for negative crossings show that in this case
the cube is not changed at all: a minus sign which appears for
$\Delta$ is then cancelled by the permutation of the two circles
(left-right) in the target space.

This completes the proof of the lemma.
\end{proof}

This lemma means that the homology of a virtual diagram with two
classical crossings (if well-defined) can be restored from an atom
endowed with orienation of the link components.

Thus, to prove that the cube of a diagram $K'$ anticommutes, we
can make some preliminary virtualizations for vertices of $K$.

To check the anticommutativity of the cube $C(K)$, we have to
consider all $2$-faces of it.

\begin{figure}
\centering\includegraphics[width=300pt]{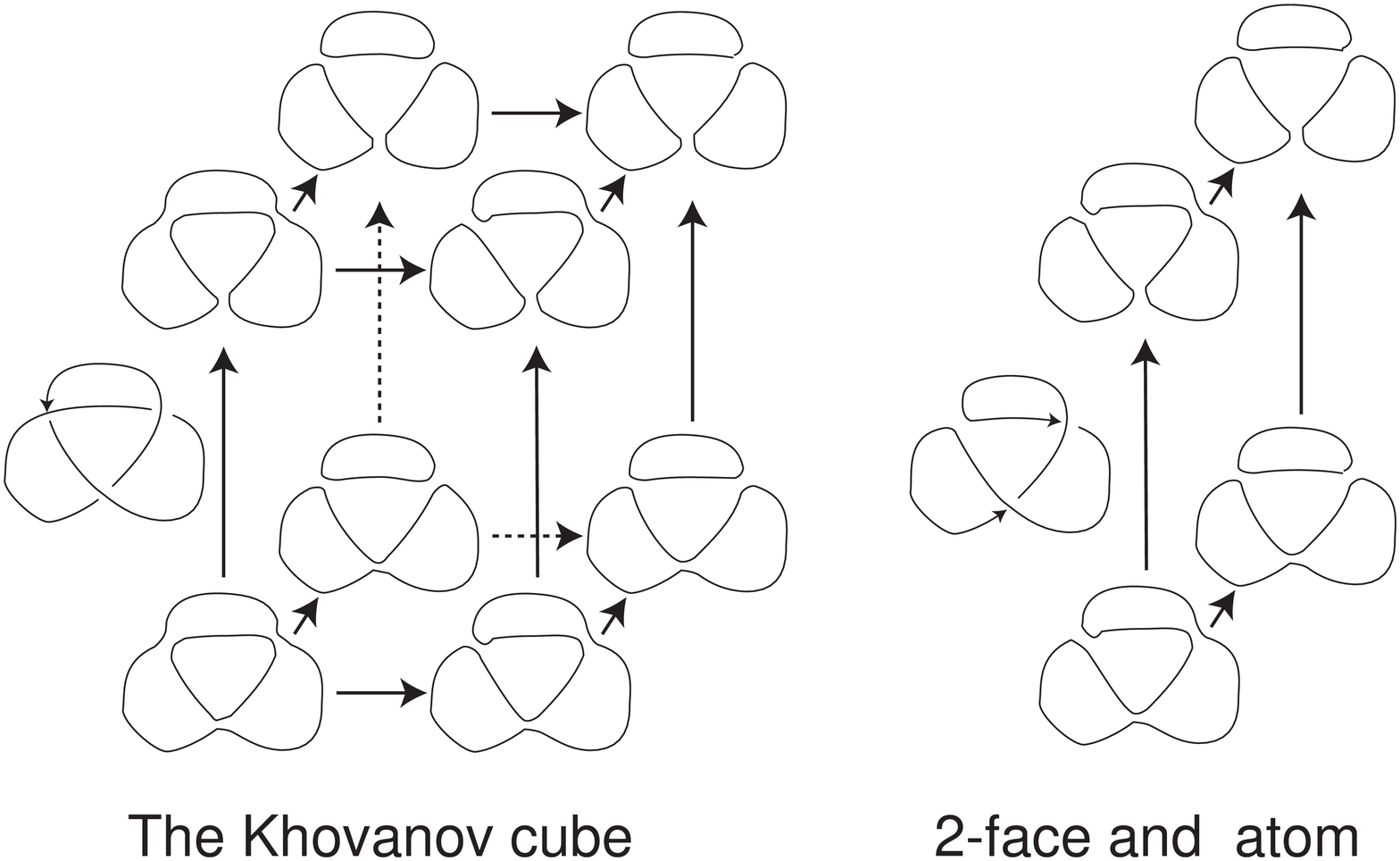} \caption{A
$2$-face generates an atom} \label{trefle}
\end{figure}

Each $2$-face is represented by fixing a way of smoothing some
$(n-2)$ vertices of $K$, see Fig.\ref{trefle}. The remaining two
crossings can be smoothed arbitrarily; the four possibilities
correspond to the vertices of the state.

Now, for these four states, there are some ``common'' circles
which do not touch any of the two vertices in question. After
removing these circles, we get an atom with two vertices.

What we actually have to check is that any face corresponding to
any possible atom with 2 vertices anti-commutes.

For two vertices of such an atom, we have some local orientations
of the link at each of these vertices; they are needed to fix the
local ordering of components while defining the differentials.

\begin{figure}
\centering\includegraphics[width=300pt]{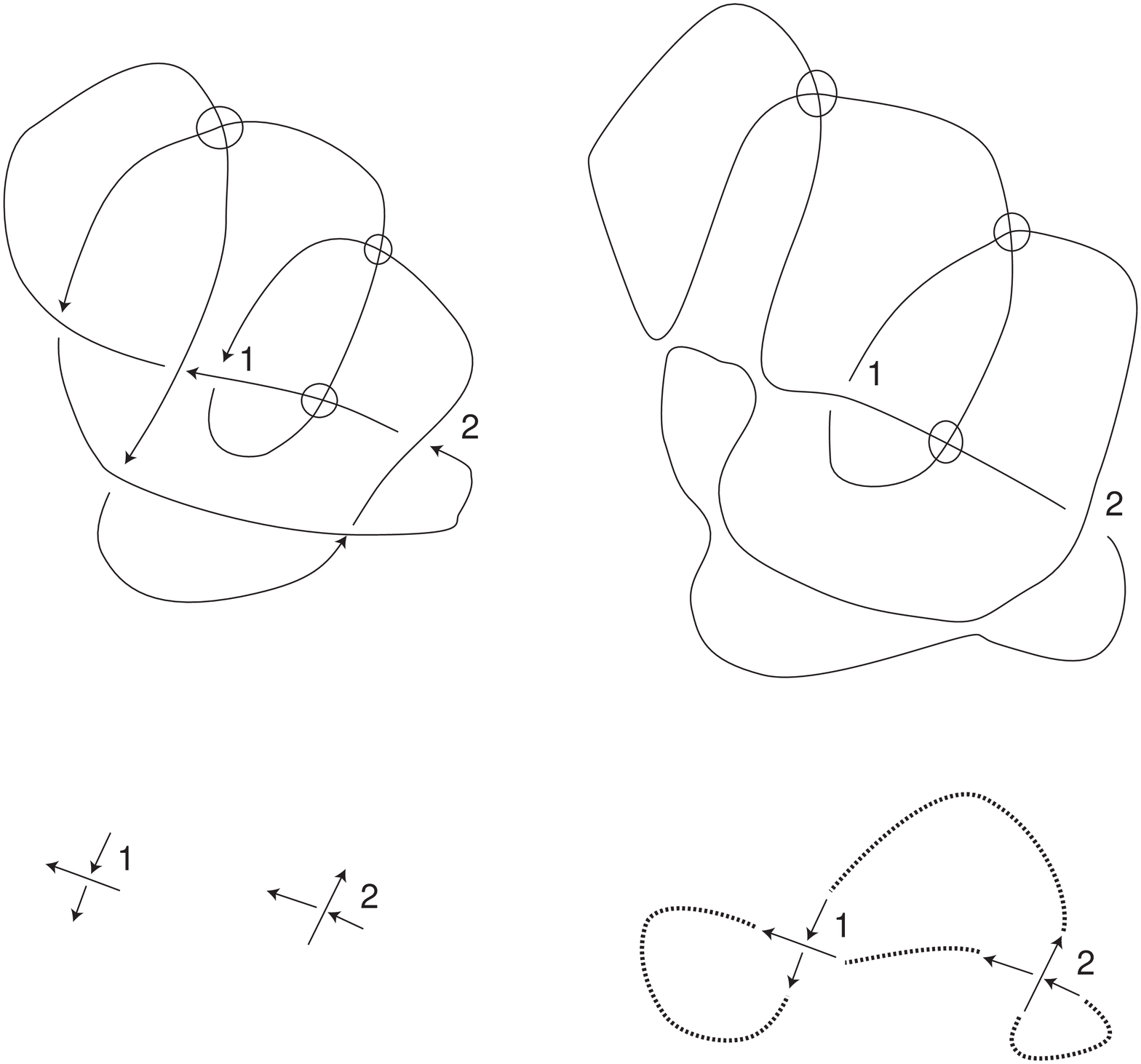}
\caption{Orientation for atom crossings} \label{opor}
\end{figure}

Note that globally these orientations might not agree on the
circles; namely, an edge of the atom with 2 vertices consists of
several edges of the diagram which might have opposite
orientations, see Fig. \ref{opor}.

It turns out, however, that these local orientations can be chosen
arbitrarily without loosing the anticommutativity property and
without changing the homology.

Fix an atom with two vertices. All possible occurences of this
atom in a cube correspond to local orientations of edges at these
vertices. Fix an orientation for one vertex $V_{1}$ and choose two
orientations for the second vertex $V_{2}$ that differ from each
other by a $\frac{\pi}{2}$-turn. Thus, we get two pictures and two
$2$-dimensional discrete cubes, $Q_1$ and $Q_2$.

\begin{lm}
If $Q_1$ is anticommutative, then so is $Q_2$.\label{lem2}
\end{lm}

\begin{proof}
Indeed, after a clockwise rotation at $V_{2}$, we change the sign
of $X$ at this crossing [again, we consider two complexes and
identify their chains so that all differentials corresponding to
the remaining vertices coincide]. This would lead to the overall
sign change of $\Delta$ if we dealt with the the usual tensor
powers.

Now, for a positive crossing, the sign of all $m$'s corresponding
to it changes as well.

For a negative crossing, the sign of all multiplications $(m)$ is
preserved whence the sign of $\Delta$ is changed for the second
time. So, the situation is as in Lemma \ref{lem1}.
\end{proof}

This lemma means that in order to check the anticommutativity of
all possible faces, it is sufficient to enumerate all atoms with 2
vertices and check the commutativity for each of them: we first
fix an representation of such an atom in ${\bf R}^{2}$ [i.e., an
immersion of its frame preserving the $A$-structure]; such
immersions differ by a possible virutalization (which does change
the homology by Lemma \ref{lem1}); then we choose a local
orientation, which does not matter either by Lemma \ref{lem2}.

First, note that among atoms with two vertices, there are
disconnected ones, i.e. those for which no edge connects one
vertex to the other.

For such atoms, if we dealt with usual tensor powers of $V$, we
would evidently get commutative faces. But the ordered tensor
powers obviously make the situation anticommute.

Some (connected) atoms with two vertices are inessential in the
following sense. We have defined the $1\to 1$ differential to be
zero. By parity reasons, in the $2$-face of an atom there might be
$0,2,$ or $4$ such edges. The case when we have no such edges is
orientable. When we have $4$ all-to-zero mappings, there is
nothing to prove. There is nothing to prove, either, when there is
one all-to-zero mapping in one composition and one all-to-zero
mapping in the other composition.

\begin{figure}
\centering\includegraphics[width=300pt]{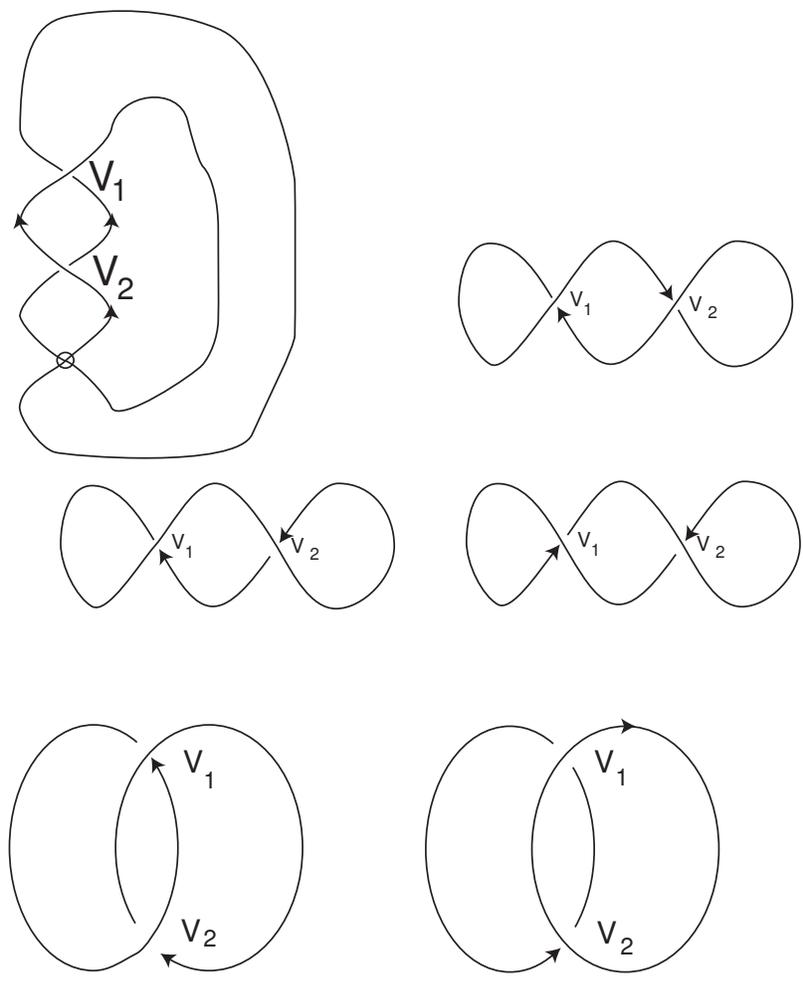}
\caption{Essential atoms with 2 vertices} \label{atoms}
\end{figure}

There are six essential connected atoms with two vertices, as
shown in Fig. \ref{atoms}. All atoms except the first one are
orientable.

\begin{figure}
\centering\includegraphics[width=300pt]{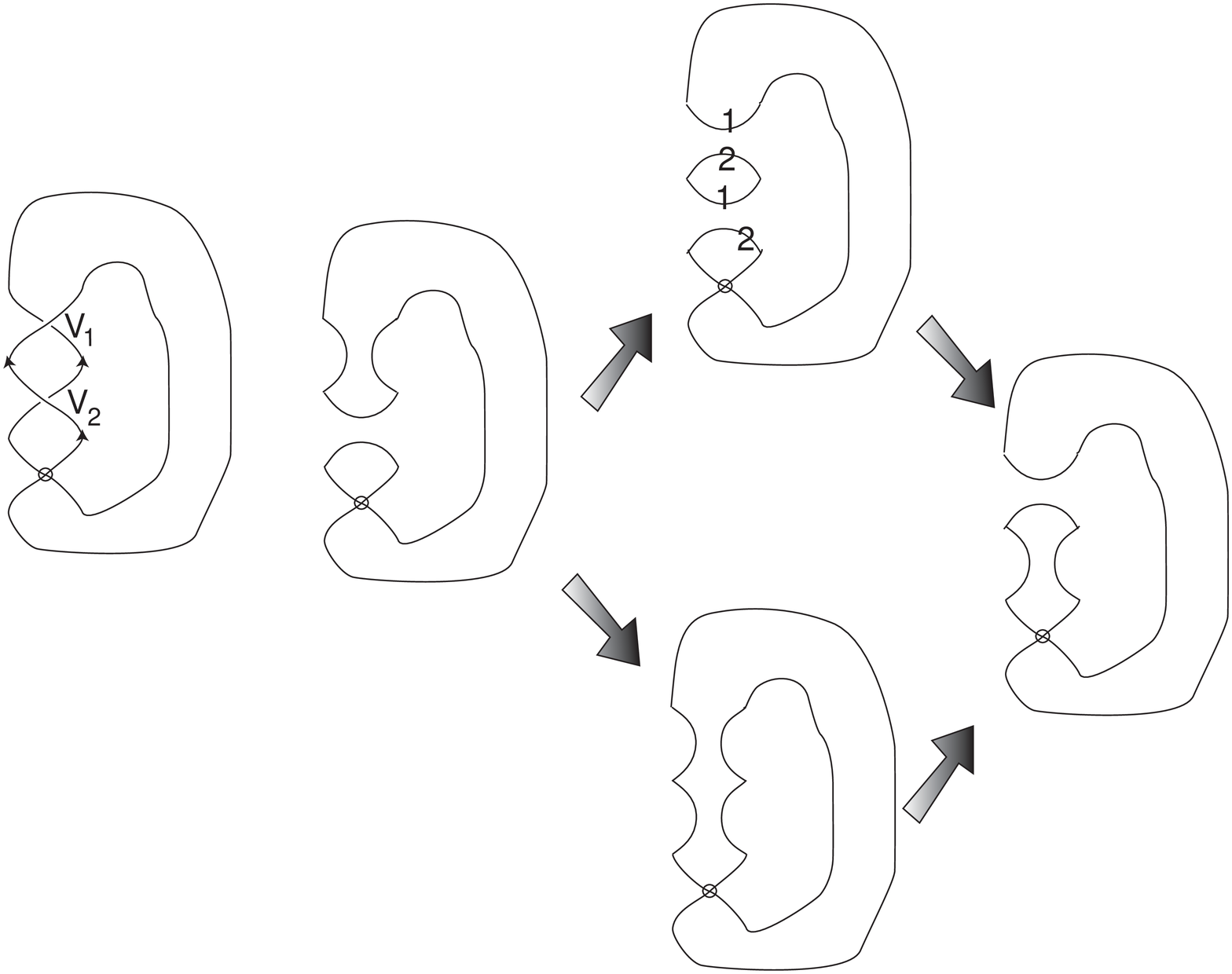} \caption{The
non-orientable atom} \label{norient}
\end{figure}

For the first one, an accurate calculation corresponding Fig.
\ref{norient} shows that both compositions give zero.

Indeed, starting with $X$, we get $\pm X\wedge X$ at the first
step and $0$ at the second step. If we start with $1$, we get
$1_{1,V_1}\wedge X_{2,V_2}+X_{1,V_1}\wedge 1_{2,V_2}$; here the
first index is the number of the circle, and the second index is
the name of the vertex. For the second vertex $V_{2}$, the first
and the second circle change their roles: $1$ becomes the lower
one and $2$ becomes the upper one. Also, for the second circle,
$X$ changes to $-X$. Thus we get $X\wedge 1- 1\wedge X$ which is
transformed by $m$ to zero.

For all orientable atoms, we fix the orientation  as shown in
Fig.\ref{atoms}. Then, the the bases  $\{1,X\}$ for vertices taken
according to \ref{figj}, will define an orientation of any circle.

Now, the anticommutativity is checked as follows. If we dealt with
the tensor product case, everything would commute. Now, the
enumeration of circles might cause minus signs on some edges. We
have to check that for any of these five atoms the total sign
would be minus.

\begin{figure}
\centering\includegraphics[width=300pt]{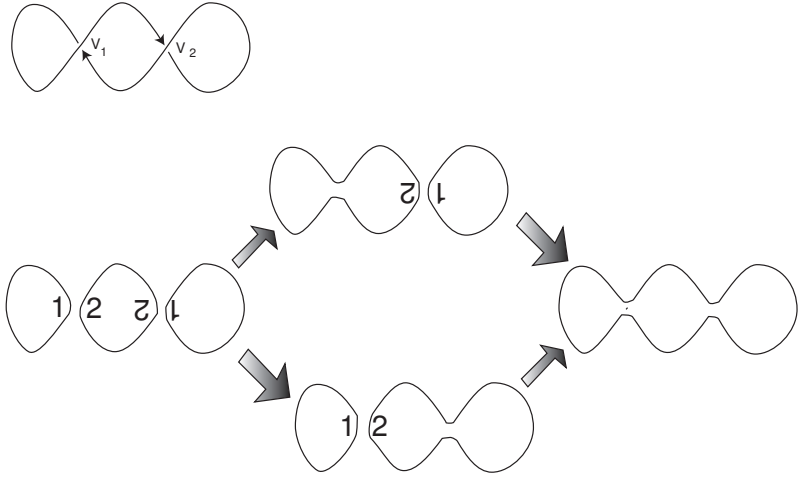} \caption{An
orientable $2$-vertex atom} \label{atom1}
\end{figure}

For instance,  in Fig. \ref{atom1} we have an oriented atom with
two vertices. The analogous check of the unordered tensor product
case means the usual associativity $m\circ(m\otimes 1)=m\circ
(1\otimes m)$, where the circles are numbered from the left to the
right.

In the left part of the figure, one pair 1\;2 is drawn upside down
to underline which circle is assumed to be locally the first
(left); the other one is the second (right).

Here we have to take into account the global ordering of the
components. Note that for three components, we have to apply
always $m\wedge Id$ first, taking those components to be
multiplied to the first and the second position.

Thus, $m\circ (m\wedge Id)$ applied to $A_{1}\wedge A_{2}\wedge
A_{3}$ gives us $m(m(A_{1},A_{2}),A_{3})=-(A_{1}\cdot A_{2}\cdot
A_{3})$; here $\cdot$ means the usual multiplication in Khovanov's
sense: $X\cdot X=0;X\cdot 1=1\cdot X=X;1\cdot 1=1$. Here the minus
sign appears because at the second crossing, we have two branches
oriented downwards, thus, the rightmost circle occurs to be
locally the left one.

On the other hand, if we consider the second crossing $V_{2}$
first, we get $A_{1}\wedge A_{2}\wedge A_{3}=(A_{2}\wedge
A_{3})\wedge A_{1}=-(A_{3}\wedge A_{2})\wedge A_{1}\to
-(A_{2}\cdot A_{3})\wedge A_{1}=A_{1}\wedge (A_{2}\cdot A_{3})$.
Applying $m$ to that, we get $A_{1}\cdot A_{2}\cdot A_{3}$.

All other atoms are checked analogously. Note that our setup gives
directly an anticommutative cube, unlike the Khovanov original
setup, where we got an anticommutative cube from a commutative one
by adding some minus signs on vertices.

\begin{thm}
$Kh(K)$ is invariant under Reidemeister moves.
\end{thm}

In fact, this theorem goes in the same lines as in \cite{BN}; one
should just take care about signs of some local differentials.

\begin{thm}
Let $K$ be a virtual diagram with orientable corresponding atom.
Then the homology $Kh(K)$ coincides with the usual Khovanov
homology \cite{Kho}.
\end{thm}

During the proof of this theorem, we denote our complex and our
homology by $C(K)$ and $Kh(K)$, and the ones constructed in
\cite{Kho} by $C'(K)$ and $Kh'(K)$, respectively.

\begin{proof}
First, we assume the diagram of $K$ is chosen in such a way that
all $X$'s for all vertices and circles agree. This is possible
since the atom corresponding to $K$ is orientable.

After that, we should just care about signs of local differential
and enumeration of circles for any vertex.

The plan is to construct a homology-preserving mapping between the
two cubes. First, $C'(K)$ does depend on enumeration of crossings.
Let us fix such an enumeration. With this enumeration, associate a
maximal spanning tree for the cubes $C(K)$ and $C'(K)$. This tree
consists of all edges of type
$(\alpha_{1},\dots,\alpha_{k},*,0,\dots,0), \alpha_{j}\in
\{0,1\}$, that is, an edge in direction $x_{j}$ is chosen if and
only if all other coordinates $x_{j}+1,\dots, x_{n}$ are zero, see
Fig.\ref{kub}.

\begin{figure}
\centering\includegraphics[width=250pt]{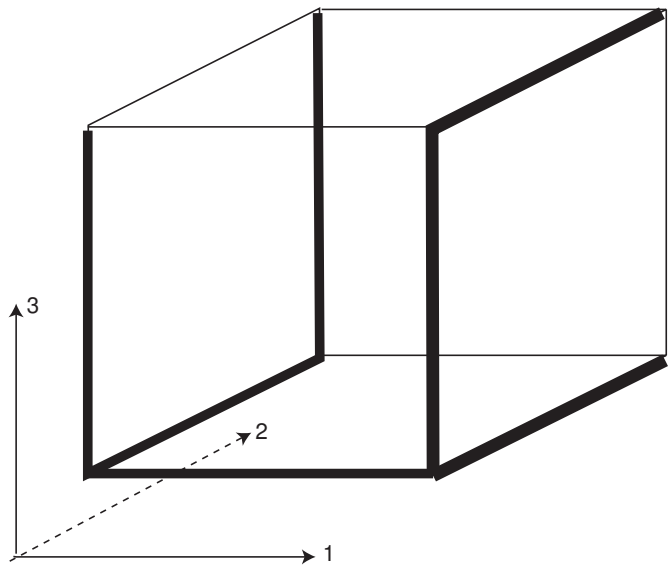} \caption{Choosing
a spanning tree} \label{kub}
\end{figure}

Now, for any state $s$ of the cube $C(K)$, we have the signed
tensor power $V^{\wedge k}$ and $V^{\otimes k}$ for $C'(K)$, where
$k$ denotes the number of circles. Enumerate the circles in the
$A$-state (the lowest corner of the cube) somehow. This ordering
defines a mapping from the space corresponding to the $A$-state
$s$ in $C(K)$ to the space corresponding to $C'(K)$. We may
prolongate this mapping to all states of the cube so that it
agrees with the differentials local along the edges of the
spanning tree.

Now, this mapping indeed agrees with all the remaining edges
because of the anticommutativity of both cubes.

\end{proof}

\section{Post Scriptum}

The construction presented above gives a partial solution to the
question: Is it true that any two classical link diagrams $K$ and
$K'$, which can be connected by a sequence of generalized
Reidemeister moves and virtualizations, are classically
equivalent. Indeed, such diagrams should have the same Khovanov
homology.

Several constructions and results concerning the Khovanov homology
generalize straightforwardly for the theory presented here.

For instance, so is the spanning tree expansion, see \cite{Weh}
and all minimality results, see \cite{Dkld}.

More precisely, the thickness of the Khovanov homology of a
virtual knot does not exceed $2+g$, where $g$ is the genus of any
atom representing $K$.

\newcommand{\RR}{\cal R}
\renewcommand{\AA}{\cal A}
\newcommand{\grad}{\mbox{deg}}

Also, Khovanov's Frobenius theory with basic ring $\RR$ and the
homology of the unknot $\AA$,

\begin{enumerate}

\item $\RR={\bf Z}[h,t]$.

\item $\AA=\RR[X]\slash (X^{2}-hX-t),$

\item $\grad X=2,\grad h=2, \grad t=4$;

\item $\Delta(1)=1\otimes X+X\otimes 1 - h 1\otimes 1$

\item $\Delta(X)=X\otimes X+t 1\otimes 1$,

\end{enumerate}

admits a straightforward generalization as above for $h=1$ [with
$1\to 1$-differential being zero. In particular, it leads to a
generalization of Lee's theory \cite{Lee}.

We shall discuss other aspects of this theory, in particular, the
general case of Frobenius' extensions and the relation of this
theory to Bar-Natan's theory for tangles and cobordisms, see
\cite{BN2,TuTu}, and other connections to Khovanov-Rozansky
homology theory.

I am very grateful to Oleg Viro for many encouraging discussions.
Also, I express my gratitude to Louis Kauffman and Victor
Vassiliev for fruitful consultations.

\end{document}